\documentclass[a4paper]{article}
\usepackage{amsmath,amssymb,mathrsfs,graphicx,psfrag}
\DeclareMathOperator{\re}{Re}
\DeclareMathOperator{\im}{Im}

\newcommand{\memo}[1]{}
\begin{document}
\title{A numerical method of computing oscillatory integral related to
hyperfunction theory}
\author{Hidenori Ogata%
\footnote{
Department of Computer and Network Engineering, Graduate School of
Informatics and Engineering, The University of Electro-Communication, 
1-5-1 Chofugaoka, Chofu, Tokyo 182-8585, Japan. 
(e-mail) {\tt ogata@im.uec.ac.jp}
}}
\maketitle
\begin{abstract}
 In this paper, we propose a numerical method of computing an integral
 whose integrand is a slowly decaying oscillatory function. 
 In the proposed method, we consider a complex analytic function in the
 upper-half complex plane, which is defined by an integral of the
 Fourier-Laplace transform type, and we obtain the desired integral by
 the analytic continuation of this analytic function onto the real axis
 using a continued fraction. 
 We also remark that the proposed method is related to hyperfunction
 theory, a theory of generalized functions based on complex function
 theory. 
 Numerical examples show the effectiveness of the proposed method. 
\end{abstract}
\section{Introduction}
\label{sec:introduction}
In this paper, we consider the computation of an integral
\begin{equation}
 \label{eq:oscillatory-integral}
  I = \int_0^{\infty}f(x)\mathrm{d}x, 
\end{equation}
where $f(x)$ is a slowly decaying oscillatory function, such as
\begin{equation*}
 \int_0^{\infty}\frac{x\cos x}{1 + x^2}\mathrm{d}x, \quad 
  \int_0^{\infty}J_{\nu}(x)\mathrm{d}x.
\end{equation*}
It is difficult to compute this integral by conventional methods such as
the DE formula \cite{TakahasiMori1974}. 
In this paper, we propose a numerical method for computing integrals of
this type. 
In the proposed method, we consider an integral of the Fourier-Laplace
transform type
\begin{equation}
 \label{eq:defining-function}
  \mathscr{F}[f](\zeta) = 
  \int_0^{\infty}f(x)\mathrm{e}^{\mathrm{i}\zeta x}\mathrm{d}x, 
\end{equation}
which is a complex analytic function in the upper-half complex plane
$\im\zeta>0$ under some conditions on $f(x)$. 
We obtain the function $\mathscr{F}[f](\zeta)$ in the upper-half
plane, and we compute the analytic continuation of
$\mathscr{F}[f](\zeta)$ on the real axis $\mathbb{R}$. 
Then, we obtain the desired integral by 
\begin{equation}
 \label{eq:boundary-value}
 I = \int_0^{\infty}f(x)\mathrm{d}x = 
  \lim_{\epsilon\rightarrow 0+}\mathscr{F}[f](\mathrm{i}\epsilon). 
\end{equation}

Our method is related to hyperfunction theory, 
a theory of generalized functions proposed by Sato \cite{Sato1959}. 
Roughly speaking, a hyperfunction is given by the boundary values of functions which are analytic 
in subdomains of $\mathbb{C}\setminus\mathbb{R}$. 
From (\ref{eq:boundary-value}), the desired integral $I$ is a boundary value on $\mathbb{R}$ 
of the function $\mathscr{F}[f](\zeta)$ analytic in the upper half plane $\im\zeta>0$, and, in this sense, 
the integral $I$ can be regarded as a value of a hyperfunction. 

Previous studies related to this paper are as follows. 
Toda and Ono proposed a method of computing
(\ref{eq:oscillatory-integral}) by the limit
\begin{equation*}
 I = 
  \lim_{n\rightarrow\infty}\int_0^{\infty}
  f(x)\exp(-2^{-n}x)\mathrm{d}x, 
\end{equation*}
which is obtained by the DE formula and the Richardson extrapolation \cite{TodaOno1978}. 
Sugihara improved Toda and Ono's method and proposed a method of
computing (\ref{eq:oscillatory-integral}) by the limit
\begin{equation*}
 I = 
  \lim_{n\rightarrow\infty}\int_0^{\infty}
  f(x)\exp(-2^{-n}x^2)\mathrm{d}x, 
\end{equation*}
which is obtained by the DE formula and the Richardson extrapolation. 
Ooura and Mori proposed a DE-type formula designed especially for
Fourier transform type integrals using a unique DE transform
\cite{OouraMori1991} and improved it for integrands with singularities
in the complex plane \cite{OouraMori1999}. 
The author proposed a method of computing Fourier transform in a way
similar to that of this paper \cite{Ogata2018}.

The remainder of this paper is structured as follows. 
In Section \ref{sec:method}, we propose a numerical method for computing
the oscillatory integral (\ref{eq:oscillatory-integral}). 
In Section \ref{sec:hyperfunction}, we give a brief review of
hyperfunction theory and mention a relation between the proposed method
and hyperfunction theory. 
In Section \ref{sec:summary}, we give the summary and concluding remarks of this paper 
and refer to problems for future study. 
\section{Numerical method}
\label{sec:method}
We consider an integral (\ref{eq:oscillatory-integral}), where $f(x)$ is
a slowly decaying oscillatory function and propose a numerical method
for computing it. 
As mentioned in the previous section, 
\begin{enumerate}
 \item we obtain a function
       $\mathscr{F}[f](\zeta)$, which is a complex analytic function in the
       upper half plane $\im\zeta>0$ under some conditions on $f(x)$, 
 \item and we obtain the desired integral $I$ by the analytic continuation of
       $\mathscr{F}[f](\zeta)$ onto the real axis. 
\end{enumerate}

The details of each step are as follows. 
\begin{enumerate}
 \item (Computation of $\mathscr{F}[f](\zeta)$ in the upper half plane
       $\im\zeta>0$)
       We obtain the analytic function $\mathscr{F}[f](\zeta)$ as a
       Taylor series
       \begin{equation*}
	\mathscr{F}[f](\zeta) = \sum_{n=0}^{\infty}c_n(\zeta -
	 \zeta_0)^n, 
       \end{equation*}
       where $\zeta_0$ is a complex constant given by the user such that 
       $\im\zeta_0 > 0$. 
       The Taylor coefficient $c_n$ is given by the integral
       \begin{equation}
	\label{eq:Taylor-coefficient}
	 c_n = 
	 \frac{1}{n!}
	 \left.\frac{\mathrm{d}^n}{\mathrm{d}\zeta^n}
	 \mathscr{F}[f](\zeta)\right|_{\zeta=\zeta_0}
	 = 
	 \frac{1}{n!}\int_0^{\infty}
	 (\mathrm{i}x)^n f(x)\mathrm{e}^{\mathrm{i}\zeta_0 x}
	 \mathrm{d}x, 
       \end{equation}
       This integral involves an exponentially decaying factor 
       $\mathrm{e}^{\mathrm{i}\zeta_0 x}$, where we remark that
       $\im\zeta_0>0$, and it is easy to compute it by conventional
       numerical integration formulas. 
 \item (Analytic continuation of $\mathscr{F}[f](\zeta)$) 
       We obtain the analytic continuation of $\mathscr{F}[f](\zeta)$
       onto the real axis $\mathbb{R}$ using its continued fraction
       expansion since the continued fraction expansion of an analytic
       function has a larger convergence region compared with that of
       its Taylor series. 
       We can convert a Taylor series 
       \begin{equation*}
	\mathscr{F}[f](\zeta) = \sum_{n=0}^{\infty}c_n(\zeta - \zeta_0)^n
       \end{equation*}
       into a continued fraction
       \begin{equation}
	\label{eq:continued-fraction}
	 \mathscr{F}[f](\zeta) = 
	 \cfrac{a_0}{1 + 
	 \cfrac{a_1(\zeta - \zeta_0)}{1 + 
	 \cfrac{a_2(\zeta - \zeta_0)}{1 + \ddots}}}
       \end{equation}
       using the so-called quotient difference (QD) algorithm (see
       Theorem 12.4c of \cite{Henrici1977}). Namely, we generate the
       sequences $e_k^{(n)}$ $( \: n = 0, 1, 2, \ldots; \: k = 0, 1, 2,
       \ldots \: )$ and $q_k^{(n)}$ $( \: n = 0, 1, 2, \ldots; \: k = 1, 2,
       \ldots \: )$ by the recurrence relations 
       \begin{gather*}
	e_0^{(n)} = 0, \quad 
	q_1^{(n)} = \frac{c_{n+1}}{c_n} \quad ( \: n = 0, 1, 2, \ldots
	\: ), 
	\\ 
	\begin{split}
	 e_k^{(n)} = q_{k}^{(n+1)} - q_{k}^{(n)} + e_{k-1}^{(n+1)},
	 \quad 
	 q_{k+1}^{(n)} = \frac{e_k^{(n+1)}}{e_k^{(n)}}q_k^{(n+1)} 
	 \\ 
	 ( \: n = 0, 1, 2, \ldots; \: k = 1, 2, \ldots \: ),
	\end{split}
       \end{gather*}
       and, then, we obtain
       \begin{equation*}
	\mathscr{F}[f](\zeta) = 
	 \cfrac{c_0}{1 - 
	 \cfrac{q_1^{(0)}(\zeta-\zeta_0)}{1 - 
	 \cfrac{e_1^{(0)}(\zeta-\zeta_0)}{1 - 
	 \cfrac{q_2^{(0)}(\zeta-\zeta_0)}{1 - 
	 \cfrac{e_2^{(0)}(\zeta-\zeta_0)}{1 - \ddots}}}}}.
       \end{equation*}
       Once we obtain the coefficients $a_n$ of the continued fraction
       (\ref{eq:continued-fraction}), we can compute the continued
       fraction by 
       \begin{equation*}
	\mathscr{F}[f](\zeta) = 
	 \lim_{k\rightarrow\infty}\frac{P_k(\zeta)}{Q_k(\zeta)},
       \end{equation*}
       where $P_k(\zeta)$ and $Q_k(\zeta)$ are the polynomials given by
       the recurrence relations
       \begin{gather*}
	P_{-1} = 0, \quad Q_{-1} = 1, \quad P_0 = c_0, \quad Q_0 = 1, 
	\\ 
	\begin{cases}
	 P_k(\zeta) = a_k(\zeta-\zeta_0)P_{k-2}(\zeta) + P_{k-1}(\zeta) 
	 \\ 
	 Q_k(\zeta) = a_k(\zeta-\zeta_0)Q_{k-2}(\zeta) + Q_{k-1}(\zeta)
	\end{cases}
	\quad ( \: k = 1, 2, \ldots \: ).
       \end{gather*}
       In the procedure of the QD algorithm, we use multiple precision
       arithmetic since the algorithm is numerically unstable. 
\end{enumerate}
Finally, we obtain the desired oscillatory integral by 
\begin{equation}
 \label{eq:boundary-value2}
  I = \int_0^{\infty}f(x)\mathrm{d}x = 
  \mathscr{F}[f](0) = 
  \lim_{\epsilon\rightarrow 0+}\mathscr{F}[f](\mathrm{i}\epsilon).
\end{equation}
\section{Relation with hyperfunction theory}
\label{sec:hyperfunction}
In this section, we give a brief review of hyperfunction theory and 
its relation with the proposed method. 
For the detail of hyperfunction theory, see the textbook \cite{Graf2010}. 

\begin{figure}[htbp]
 \begin{center}
  \begin{tabular}{cc}
   \psfrag{x}{$\re z$}
   \psfrag{y}{$\im z$}
   \psfrag{z}{\rotatebox{90}{\hspace{-5mm}$\re F(z)$}}
   \includegraphics[width=0.49\textwidth]{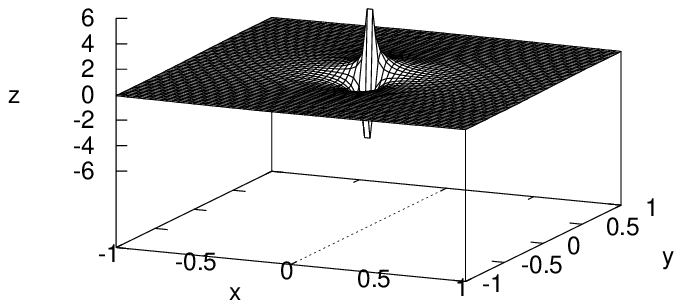}& 
       \psfrag{x}{$\re z$}
       \psfrag{y}{$\im z$}
       \psfrag{z}{\rotatebox{90}{\hspace{-5mm}$\re F(z)$}}
       \includegraphics[width=0.49\textwidth]{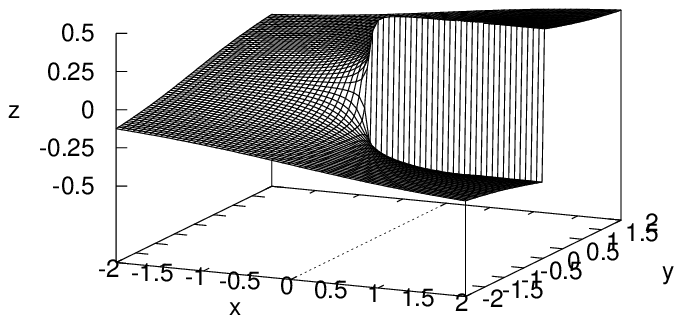}
       \\
   $\delta(x)$ & $Y(x)$
  \end{tabular}
 \end{center}
 \caption{The graphs of the defining functions $F(z)$ of the Dirac delta function $\delta(x)$ 
 and the Heaviside step function $Y(x)$.}
 \label{fig:defining-function}
\end{figure}

Intuitively, a hyperfunction is defined by the difference between the
boundary values of a complex analytic functions $F_+(z)$ in
a subdomain of the upper half plane $\im z>0$ and a complex analytic function 
$F_-(z)$ in a subdomain of the lower half plane $\im z<0$ 
\begin{equation*}
 \mathfrak{f}(x) = F_+(x+\mathrm{i}0) - F_-(x-\mathrm{i}0), 
\end{equation*}
where the analytic functions $F_{\pm}(z)$ are called defining functions
of the hyperfunction $\mathfrak{f}(x)$. 
We also use the notation 
\begin{equation*}
 \mathfrak{f}(x) = [F(z)] = F(x+\mathrm{i}0) - F(x-\mathrm{i}0), 
\end{equation*}
where 
\begin{equation*}
 F(z) = 
  \begin{cases}
   F_+(z) & \im z > 0 \\ 
   F_-(z) & \im z < 0.
  \end{cases}
\end{equation*}
For example, the Dirac delta function $\delta(x)$ is defined as a
hyperfunction by 
\begin{equation*}
 \delta(x) = 
  - \frac{1}{2\pi\mathrm{i}}
  \left( \frac{1}{x+\mathrm{i}0} - \frac{1}{x-\mathrm{i}0} \right), 
\end{equation*}
and the Heaviside step function $Y(x)$ is defined as a hyperfunction by 
\begin{equation*}
 Y(x) = 
  - \frac{1}{2\pi\mathrm{i}}
  \left\{\log(-(x+\mathrm{i}0)) - \log(-(x-\mathrm{i}0)) \right\} 
\end{equation*}
where the complex logarithmic function $\log z$ is the principal value,
that is, the branch such that it takes a real value on the positive real
axis. 
Figure \ref{fig:defining-function} shows the graphs of the defining functions 
$F(z)$ of the Dirac delta function $\delta(x)$ and the Heaviside step function $Y(x)$. 
From these figures, we find easily that the difference between the boundary values on the real axis 
of the defining functions $F(z)$ gives the hyperfunctions. 

More precisely, hyperfunctions are defined as follows. 
Let $I$ be an open interval on the real axis $\mathbb{R}$ and 
$D$ be a complex neighborhood of $I$, that is, a complex domain including $I$ as a closed subset. 
We denote the set of complex analytic functions in $D\setminus I$ by $\mathscr{O}(D\setminus I)$, 
which can be regarded as a complex linear space by the usual summation and multiplication by complex numbers, 
and the set of complex analytic functions in $D$ by $\mathscr{O}(D)$, 
which can be regarded as a linear subspace of $\mathscr{O}(D\setminus I)$. 
Then, we can consider the quotient space 
$\mathscr{B}(I)\equiv\mathscr{O}(D\setminus I)/\mathscr{O}(D)$. 
We call an element of $\mathscr{B}(I)$, that is, an equivalence class $[F(z)]$ 
$( \: F(z) \in \mathscr{O}(D\setminus I) \: )$ a hyperfunction on $I$ and denote it by 
\begin{equation*}
 \mathfrak{f}(x) = [F(z)] = F(x+\mathrm{i}0) - F(x-\mathrm{i}0), 
\end{equation*}
where we call $F(z)$ a defining function of the hyperfunction $\mathfrak{f}(x)$. 
We also use the notation
\begin{equation*}
 \mathfrak{f}(x) = [F_+(z), F_-(z)] = F_+(x+\mathrm{i}0) - F_-(x-\mathrm{i}0), 
\end{equation*}
where
\begin{equation*}
 F_{\pm}(z) = F|_{D_{\pm}}(z) \quad \mbox{with} \quad 
  D_{\pm} = \{ \: z \in \mathbb{C} \: | \: \pm\im z > 0 \: \}. 
\end{equation*}
In other words, a hyperfunction $\mathfrak{f}(x)$ on $I$ is a complex analytic function $F(z)$ in 
$D\setminus I$, where another complex analytic function $\widetilde{F}(z)$ in $D\setminus I$ is 
identified with $F(z)$ if $\widetilde{F}(z) - F(z)$ is analytic in $D$. 

We define the value of a hyperfunction $\mathfrak{f}(x)=F_+(x+\mathrm{i}0) - F_-(x-\mathrm{i}0)$ at a point 
$x_0\in I$ by 
\begin{equation*}
 \mathfrak{f}(x_0) = 
  \lim_{\epsilon\rightarrow 0+}
  \left\{ F_+(x_0+\mathrm{i}\epsilon) - F_-(x_0-\mathrm{i}\epsilon)\right\}
\end{equation*}
if the limit on the right hand side exists. 
We do not define the value $\mathfrak{f}(x_0)$ if the limit does not exist. 

We return to the discussion on the method of computing the oscillatory integral (\ref{eq:oscillatory-integral}). 
From (\ref{eq:boundary-value2}) and the fact that 
$\mathscr{F}[f](\zeta)$ is a complex analytic function in the upper half plane 
$\im\zeta > 0$, 
we can regard the desired oscillatory integral $I$ as the value 
at $\xi=0$ of the hyperfunction 
$\mathfrak{f}(x)=[F_+(\zeta), F_-(\zeta)]$ 
whose defining functions are 
$F_+(\zeta) = \mathscr{F}[f](\zeta)$ and $F_-(\zeta) = 0$. 
\section{Numerical examples}
\label{sec:example}
In this section, we discuss show some numerical examples that illustrate 
the effectiveness of the proposed method. 
All the computations were performed using programs coded in C++ with 100
decimal digit precision, where we use the multiple precision arithmetic
library {\it exflib} \cite{Fujiwara-exflib}.

We computed the oscillatory integrals
\begin{align*}
 \mathrm{(1)} \quad & 
 \int_0^{\infty}\frac{\cos(x/2)-\cos x}{x}\mathrm{d}x = \log 2 =
 0.69315\ldots, 
 \\ 
 \mathrm{(2)} \quad & 
 \int_0^{\infty}\log x\cos x\mathrm{d}x = - \frac{\pi}{2} = 
 - 1.57080\ldots,
 \\
 \mathrm{(3)} \quad & 
 \int_0^{\infty}J_0(x)\mathrm{d}x = 1, 
 \\ 
 \mathrm{(4)} \quad & 
 \int_0^{\infty}\frac{xJ_0(x)}{x^2+1}\mathrm{d}x = K_0(1) =
 0.42102\ldots, 
 \\
 \mathrm{(5)} \quad & 
 \int_0^{\infty}\frac{J_0(x)}{(x^2+1)^{1/2}}\mathrm{d}x = 
 K_0\left(\frac{1}{2}\right)I_0\left(\frac{1}{2}\right) = 0.98310\ldots, 
 \\
 \mathrm{(6)} \quad & 
 \int_0^{\infty}\log x J_0(x)\mathrm{d}x = - \gamma - \log 2 = 
 - 1.2704\ldots, 
 \\
 \mathrm{(7)} \quad & 
 \int_0^{\infty}\frac{x J_1(\sqrt{x^2+1})}{\sqrt{x^2+1}}\mathrm{d}x =
 J_0(1) = 0.76520\ldots, 
 \\ 
 \mathrm{(8)} \quad & 
 \int_0^{\infty}\frac{Y_0(x)}{x^2+1}\mathrm{d}x = - K_0(1) = 
 - 0.42102\ldots.
\end{align*}
by the proposed method. 
We took the center $\zeta_0$ of the Taylor series 
\begin{math}
 \mathscr{F}[f](\zeta) = \sum_{n=0}^{\infty}c_n(\zeta-\zeta_0)^n
\end{math}
as $\zeta_0 = \mathrm{i}$, and we computed the Taylor coefficients 
$c_0, c_1, c_2, \ldots, c_{100}$ given by the integral
(\ref{eq:Taylor-coefficient}), 
which were computed using the DE formula. 
Table \ref{tab:example1} show the relative error and the number of
functional evaluations for $f(x)$ in (\ref{eq:oscillatory-integral}) of
the proposed method. 
Table \ref{tab:example1} show that we can compute oscillatory integrals
with high accuracy by the proposed method.
\begin{table}[htbp]
 \caption{The relative errors and the number of functional evaluations
 for $f(x)$ in (\ref{eq:oscillatory-integral}) of the proposed method
 applied to the integrals (1-8).}
 \begin{center}
  \begin{tabular}{ccccc}
   \hline
   integral & (1) & (2) & (3) & (4) \\
   \hline
   relative error & 5.4E-26 & 6.2E-35 & 3.8E-34 & 1.4E-36 \\ 
   number of functional evaluations & 917 & 964 & 957 & 927 \\ 
   \hline
   integral & (5) & (6) & (7) & (8) \\ 
   \hline
   relative error & 1.3E-35 & 3.8E-36 & 1.1E-33 & 2.1E-37 \\ 
   number of functional evaluations & 954 & 958 & 927 & 947 \\ 
   \hline
  \end{tabular}
 \end{center}
 \label{tab:example1}
\end{table}

Here, we point out an interesting fact. Since we took
$\zeta_0=\mathrm{i}$, the integral in (\ref{eq:Taylor-coefficient}) 
which gives the Taylor coefficients $c_n$ is given by
\begin{equation*}
 c_n = 
  \frac{1}{n!}\int_0^{\infty}(\mathrm{i}x)^n
  f(x)\mathrm{e}^{-x}\mathrm{d}x, 
\end{equation*}
which does not involve any oscillatory function. 
It means that we can obtain oscillatory integrals without computing
oscillatory integrals. 
In other words, we can obtain Fourier transform type integrals by
computing Laplace transform type integrals. 

As a comparison with the results by a previous method, we computed some
of the integrals by the Euler transform for computing alternating
series. 
Namely, we express the desired oscillatory integral as 
\begin{equation}
 \label{eq:alternating-series}
  \int_0^{\infty}f(x)\mathrm{d}x = 
  \sum_{k=0}^{\infty}(-1)^k
  \left|\int_{a_k}^{a_{k+1}}f(x)\mathrm{d}x\right|, 
\end{equation}
where $0=a_0 < a_1 < a_2 < \ldots$ and 
\begin{equation*}
 \int_{a_k}^{a_{k+1}}f(x)\mathrm{d}x \cdot 
  \int_{a_{k+1}}^{a_{k+2}}f(x)\mathrm{d}x < 0 
  \quad ( \: k = 0, 1, 2, \ldots \: ).
\end{equation*}
We computed each integral on the right hand side of
(\ref{eq:alternating-series}) by the 100-point Gauss-Legendre formula, 
and we obtained the desired integral by applying the Euler transform to
the alternating series (\ref{eq:alternating-series}). 
Table \ref{tab:example2} show the relative errors of the method using
the Euler transform applied to the integrals (2,3,4,5). 
The number of functional evaluations for $f(x)$ in
(\ref{eq:oscillatory-integral}) are 5000 in all the computations of the
integrals. 
Comparing the results in Table \ref{tab:example1} and \ref{tab:example2},
the proposed method is far superior to the method using the Euler
transform. 
\begin{table}[htbp]
 \caption{The relative errors of the method using the Euler transform
 applied to the integrals (2,3,4,5). The number of functional
 evaluations for $f(x)$ in (\ref{eq:oscillatory-integral}) are 5000 in
 all the computations of the integrals.}
 \begin{center}
  \begin{tabular}{ccccc}
   \hline
   integral & (2) & (3) & (4) & (5) \\
   \hline
   relative error & 6.3E-25 & 1.9E-25 & 3.1E-25 & 1.4E-25
		   \\
   \hline
  \end{tabular}
 \end{center}
 \label{tab:example2}
\end{table}
\section{Summary}
\label{sec:summary}
In this paper, we proposed a numerical method of computing an integral involving 
a slowly decaying oscillatory function. 
In the proposed method, we introduce a complex analytic function in the upper half plane 
and obtain the desired integral by the analytic continuation of this analytic function 
onto the real axis.  
The analytic function is obtained as a Taylor series by computing its coefficients, 
which are given by integrals involving exponentially decaying functions and easily evaluated 
by the conventional numerical integration formula. 
The analytic continuation is performed by converting the Taylor series into a continued fraction 
using the quotient difference (QD) algorithm. 
Then, the desired integral is obtained as the value of the analytic function at the origin. 
We also remarked that the proposed method is closely related to hyperfunction theory in the sense that 
the desired oscillatory integral can be regarded as the boundary value of the hyperfunction whose 
defining function is the above complex analytic function. 

As presented in this paper, 
hyperfunctions, that is, the boundary values of analytic functions appear in many studies of science 
and engineering. 
In the evaluation of hyperfunctions, the techniques of analytic continuation are necessary. 
Therefore, the author believes that the study and the computation of hyperfunctions together with 
analytic continuation will be important and useful in applied mathematics. 

Problems for future study related to this paper are, for example, theoretical error estimate of 
the proposed method and the reduction of the computational cost of analytic continuation, 
which needs multiple precision arithmetic in the QD algorithm applied to the conversion of a Taylor 
series into a continued fraction. 
\bibliographystyle{plain}
\bibliography{arxiv2019_1}
\end{document}